\documentclass[oneside,a4paper]{article}
\usepackage{amssymb, amsfonts, amsmath, amsthm}
\usepackage{enumerate}
\usepackage[colorlinks=true,
citecolor=black,
linkcolor=black,
anchorcolor=black,
filecolor=black,
menucolor=black,
urlcolor=black,
pdftitle={Intrinsic isometries in Euclidean space},
pdfauthor={Anton Petrunin},
]{hyperref}

\usepackage{epsfig,lpic,wrapfig}
\renewcommand*{\HyperDestNameFilter}[1]{\jobname-#1}

\def\smallskip{\vskip\smallskipamount}
\def\medskip{\vskip\medskipamount}
\def\bigskip{\vskip\bigskipamount}

\def\emptyset{\varnothing}

\def\ae{\buildrel \hbox{\tiny\it a.e.} \over{=\joinrel=}}

\def\parit#1{\medskip\noindent{\it #1}}
\def\parbf#1{\medskip\noindent{\bf #1}}
\def\qeds{\qed\par\medskip}

\newcounter{thm}[section]

\def\claim#1{\par\medskip\noindent\refstepcounter{thm}\hbox{\bf \arabic{section}.\arabic{thm}. #1.}
\it\ %\ignorespaces
}
\def\endclaim{
\par\medskip}
\newenvironment{thm}{\claim}{\endclaim}

\def\EE{\mathbb{E}}%
\def\II{\mathbb{I}}
\def\JJ{\mathbb{J}}

\def\MM{\mathbb{M}}% 
\def\RR{\mathbb{R}}%  
\def\spc#1{\mathcal{#1}}

\def\~{\tilde}
\def\eps{\varepsilon}
\def\phi{\varphi}
\def\theta{\vartheta}

\def\ge{\geqslant}
\def\le{\leqslant}

\def\i{\subset}
\def\l{\!\left}
\def\r{\right}
\def\<{\langle}
\def\>{\rangle}

\def\:{\colon}

\def\diam{\mathop{\rm diam}\nolimits}
\def\dim{\mathop{\rm dim}\nolimits}
\def\dist{\mathop{\rm pull}\nolimits}
\def\length{\mathop{\rm length}\nolimits}
\def\pack{\mathop{\rm pack}\nolimits}

\newcommand*{\z}[1]{#1\nobreak\discretionary{}%
            {\hbox{$\mathsurround=0pt #1$}}{}}
\let\oldcdot\cdot
\def\cdot{{\hskip0.5pt\z\oldcdot\hskip0.5pt}}

\begin{document}
\title{Intrinsic isometries in Euclidean space}
\author{Anton Petrunin}
\date{}
\maketitle

\begin{abstract}
I consider compact metric spaces which admit intrinsic isometries to Euclidean $d$-space.
The main result roughly states that the class of these spaces coincides with class of inverse limits of Euclidean $d$-polyhedrons. 
\end{abstract}

\section{Introduction} 

The \emph{intrinsic isometries} are defined in Section~\ref{preliminaries};
it is a variation of notion of \emph{path isometry}, that is a map that preserves the lengths of curves.
Any intrinsic isometry is a path isometry, the converse does not hold in general.

The following statement is one of the reason we prefer intrinsic isometry.

\begin{thm}{Trivial  statement}\label{top-dim}
Let $\spc{X}$ be a compact metric space which admits an intrinsic isometry to $d$-dimensional Euclidean space (further denoted by $\EE^d$). 
Then 
$\dim \spc{X}\le d$,
where $\dim$ denotes the Lebesgue's covering dimension.
\end{thm}

This statement is proved in Section~\ref{proofs}.
An analogous statement for path isometry does not hold, see Example~\ref{exam:path}.
The Hausdorff dimension cannot be bounded on the similar way.
For example, $\RR$-tree admits an intrinsic isometry to $\RR$ and it contains compact subsets of arbitrary large Hausdorff dimension.

\medskip

Here are some known results 
on length spaces which admit \emph{intrinsic isometry} to $\EE^d$.

\begin{thm}{Theorem}\label{thm:gromov}
Let $\spc{R}$ be $d$-dimensional Riemannian space 
and $f\:\spc{R}\to\EE^d$ be a short map%
\footnote{That is, a $1$-Lipschitz map.}%
.
Then given $\eps>0$, 
there is an intrinsic isometry $\imath\:\spc{R}\to\EE^d$ such that 
$$|f(x)-\imath(x)|_{\EE^d}<\eps$$
for any $x\in \spc{R}$.

In particular, any Riemannian $d$-space admits an intrinsic isometry to $\EE^d$.
\end{thm}

For path isometries, this theorem was proved by Mikhael Gromov \cite[2.4.11]{gromov-PDE},
and the same proof works for intrinsic isometries.
Applying this theorem, one can show 
that any limit of increasing sequence of Riemannian metrics on a fixed $d$-dimensional manifold admits an intrinsic isometry to $\EE^d$.
(The proof is similar to ``if''-part of the main theorem.)
In particular, any sub-Riemannian metric on $d$-dimensional manifold admits an intrinsic isometry to $\EE^d$.

\begin{thm}{Theorem}\label{PL-Nash}
Let $\spc{P}$ be a Euclidean polyhedron and $f\: \spc{P}\to\EE^d$ be a short map.
Then, given $\eps>0$, there is a piecewise linear intrinsic isometry 
$\imath\: \spc{P}\to \EE^d$ such that
$$|f(x)-\imath(x)|_{\EE^d}<\eps$$
for any $x\in \spc{P}$.
\end{thm}

\begin{thm}{Corollary}
Any $d$-dimensional Euclidean polyhedron admits a piecewise linear intrinsic isometry to $\EE^d$.
\end{thm}

The corollary was proved by Viktor Zalgaller \cite{zalgaller} for dimension $\le 4$, 
but a slight modification of his proof works in all dimensions, 
see \cite{krat}.
The 2-dimensional case of Theorem~\ref{PL-Nash} was proved by Svetlana Krat \cite{krat}.
Later, the proof was extend to all dimensions by Arseniy Akopyan \cite{akopjan}.
His proof use a \emph{piecewise linear analog of Nash--Kuiper theorem} which was proved by Ulrich Brehm \cite{brehm},
and reproved independently by  Arseniy Akopyan and Alexey Tarasov~\cite{akopjan-tarasov}.

\parbf{Iff-condition.} Now we describe the main result of the paper.

A compact metric space $\spc{X}$ is called \emph{pro-Euclidean space of rank $\le d$} 
if it can be presented as an \emph{inverse limit} $\spc{X}=\varprojlim \spc{P}_n$ (see Section~\ref{preliminaries}) of a sequence of Euclidean $d$-polyhedrons $\spc{P}_n$.

\begin{thm}{Main theorem}\label{main}
A compact metric space $\spc{X}$ admits an intrinsic isometry to $\EE^d$ 
if and only if $\spc{X}$ is a pro-Euclidean space of rank $\le d$. 
\end{thm}

The proof is straightforward. 
I like the formulation of the theorem.
It is rare when inverse limits help to solve a natural problem in metric geometry;
the only other example I know is the characterization of  homogeneous
locally compact metric spaces given by Valerii Berestovskii~\cite{berestovskii}.

Note that the statement in theorem \ref{thm:gromov} (in the compact case) 
is equivalent to the fact that any compact Riemannian 
$d$-space is a pro-Euclidean space of rank $\le d$.
The later can be obtained directly from the following exercise;
this way the main theorem provides an alternative  proof to Theorem~\ref{thm:gromov} in the compact case.

\begin{thm}{Exercise}
Show that any compact Riemannian space admits a Lipschitz approximation by Euclidean polyhedrons.
\end{thm}

\parbf{A non-example.} 
Let us remind  that \emph{Minkowski space} is finite dimensional real vector spaces with metric induced by a norm.

\begin{thm}{Proposition}\label{minkowski}
Let $\Omega$ be an open subset of  Minkowski $d$-space $\MM^d$.
Assume $\Omega$ admits an intrinsic\footnote{In fact the same is true for path isometry.} isometry  to $\EE^m$ then $d\le m$ and $\MM^d$ is isometric to~$\EE^d$.
\end{thm}

In particular, the condition in \ref{top-dim} on Lebesgue's dimension is not sufficient.

\smallskip

I'm grateful to 
Arseniy Akopyan, 
Dmitri Burago, 
Sergei Ivanov,
Alexander Lytchak
and an anonymous referee for helpful letters and discussions.

\section{Preliminaries}\label{preliminaries}

\parbf{Standard definitions.}
Given a metric space $\spc{X}$ and two points $x,x'\in \spc{X}$, we will denote by $|x-x'|=|x-x'|_{\spc{X}}$ the distance from $x$ to $x'$ in $\spc{X}$.

A \emph{length space} is a metric space such that for any two points $x,x'$ the distance $|x - x'|$ coincides with the greatest lower bound of lengths of curves connecting $x$ and $x'$.

A map $f\:\spc{X}\to \spc{Y}$ between metric spaces $\spc{X}$ and $\spc{Y}$ is called \emph{short} if for any $x,x'\in \spc{X}$ we have
$$|f(x)-f(x')|_{\spc{Y}}\le |x -x'|_{\spc{X}}.$$

A length space
$\spc{P}$ is called \emph{Euclidean $d$-polyhedron} if there is a finite triangulation of $\spc{P}$ such that each simplex is isometric to a simplex in $\EE^d$% and maximal dimension of simplex is equal $d$
.

\parbf{Inverse limit.} 
Consider an \emph{inverse system} of compact metric spaces $(\spc{X}_n)_{n=0}^\infty$ and short maps $\phi_{m,n}:\spc{X}_m\to \spc{X}_n$ for $m\ge n$;
that is,
\begin{enumerate}
\item $\phi_{m,n}\circ \phi_{k,m}=\phi_{k,n}$ for any triple $k\ge m\ge n$ and
\item for any $n$, the map $\phi_{n,n}$ is identity map of $\spc{X}_n$.
\end{enumerate}

A compact metric space $\spc{X}$ is called \emph{inverse limit of the system $(\phi_{m,n}, \spc{X}_n)$} (denoted by $\spc{X}=\varprojlim \spc{X}_n$) if its underling space consists of all sequences $x_n\in \spc{X}_n$ such that $\phi_{m,n}(x_m)=x_n$ for all $m\ge n$ and for two such sequences $(x_n)$ and $(x'_n)$ the distance is defined by 
$$|(x_n)-(x'_n)|_{\spc{X}}=\lim_{n\to\infty}|x_n- x'_n|_{\spc{X}_n}.$$

If $\spc{X}=\varprojlim \spc{X}_n$, then the map $\psi_n\:\spc{X}\to \spc{X}_n$, defined by $\psi_n\:(x_i)_{i=0}^\infty\mapsto x_n$ are called \emph{projections}.
Clearly $\psi_n=\phi_{m,n}\circ\psi_m$ for all $m\ge n$. 

\parit{Comments.}
The above definition is equivalent to the usual inverse limit in the category 
with class of objects formed by compact metric spaces and class of morphisms formed by short maps.

Note that inverse limit is not always defined, and
if defined it is the result is compact by definition.
(In principle, the category of compact metric spaces can be extended so that the limit of any inverse system is well defined.)

It is easy to see that inverse limit of length spaces is a length space.

In general, the inverse limit of a system of spaces differ from its Gromov--Hausdorff limit.
For example, consider inverse system  $\spc{X}_n=[0,1]$ with maps $\phi_{m,n}(x)\equiv 0$. 
The inverse limit of this system is isometric to one-point space, while the Gromov--Hausdorff limit is isometric to $[0,1]$.
Nevertheless, it is easy to see that if for any $\eps>0$ the images of $\phi_{m,n}$ form an $\eps$-net in $X_n$ for all sufficiently large $m$ and $n$, then $\spc{X}=\varprojlim \spc{X}_n$ is isometric to the Gromov--Hausdorff limit.

\parbf{Intrinsic isometries and pull back metrics.}
Let $\spc{X}$ and $\spc{Y}$ be metric spaces and  $f\:\spc{X}\to \spc{Y}$ be continuous map.
Given two points $x,x'\in \spc{X}$, a sequence of points $x=x_0,x_1,\dots,x_n=x'$ is called $\eps$-chain from $x$ to $x'$ if $|x_{i-1}-x_i|\le\eps$ for all $i>0$.
Set
$$\dist_{f,\eps}(x,x')
=
\inf\l\{\sum_{i=1}^n|f(x_{i-1})-f(x_{i})|_{\spc{Y}}\r\}$$
where the greatest lower bound is taken along all  $\eps$-chains $(x_i)_{i=0}^n$ from $x$ to $x'$.

Clearly $\dist_{f,\eps}$ is a pseudometric%
\footnote{That is, it satisfies triangle inequality, it is symmetric, non-negative and $\dist_{f,\eps}(x,x)=0$, but it might happen that $\dist_{f,\eps}(x,x')=0$ for $x\not=x'$.}
on $\spc{X}$ and $\dist_{f,\eps}(x,x')$ is non-increasing in $\eps$.
Thus, the following (possibly infinite) limit
$$\dist_{f}(x,x')=\lim_{\eps\to0}\dist_{f,\eps}(x,x')$$
is well defined.
The pseudometric $\dist_f\:\spc{X}\times\spc{X}\to[0,\infty]$ will be called the \emph{pull back metric} for $f$.

A map $f\:\spc{X}\to \spc{Y}$ between length spaces $\spc{X}$ and $\spc{Y}$ is an \emph{intrinsic isometry} if 
$$|x-x'|_{\spc{X}}=\dist_f(x,x')$$
for any $x,x'\in \spc{X}$.

Any intrinsic isometry is a short map.
Moreover, it is easy to see that intrinsic isometry preserves the lengths of curves.
The converse does not hold, see Section~\ref{path.isometry}.

\begin{thm}{Proposition}\label{int>int}
Let $\spc{X}$ be a compact (or even proper\footnote{That is, all closed bounded sets in $\spc{X}$ are compact.}) metric space.
Then existence of an intrinsic isometry $f\:\spc{X}\to \spc{Y}$ implies that $\spc{X}$ is a length space.
\end{thm}

The proof is left to the reader.
Note that the proposition does not hold for general complete space $\spc{X}$.
Consider two points which connected by countable number of unit intervals $\II_n$ 
and one interval of length $\frac12$; 
equip the obtained space with the natural intrinsic metric.
Let us remove from our space the interval of length $\frac12$.
The metric on the remaining space $\spc{X}$ is not intrinsic, but complete.
Further let us construct a map $f\:\spc{X}\to\RR$ so that the restriction $f_n=f|_{\II_n}$ 
is an intrinsic isometry, $f_n(0)=0$, $f_n(1)=\tfrac12$ and $f_n(x)$ converges uniformly to $\tfrac x2$.
It is easy to see that $f\:\spc{X}\to\RR$ is an intrinsic isometry.

The following statement is analogous to the semi-continuity of length functional.

\begin{thm}{Proposition}\label{prop:pull-back} Let $\spc{X}$ and $\spc{Y}$ be metric spaces, $\spc{X}$ be compact and the continuous map $f\:\spc{X}\to \spc{Y}$ is such that
$$\sup_{x,x'\in\spc{X}}\dist_f(x,x')<\infty.$$
Then given $\eps>0$ there is $\delta=\delta(f,\eps)>0$ such that for any short map $h\:\spc{X}\to \spc{Y}$ such that
$$|f(x)-h(x)|_{\spc{Y}}<\delta\ \ \text{for any}\ \ x\in \spc{X}$$
we have
$$\dist_{f}(x,x')<\dist_{h}(x,x')+\eps$$
for any $x,x'\in \spc{X}$.
\end{thm}

The proof is a direct application of the Lemma~\ref{lem:pull-back}.

For a compact metric space $\spc{X}$, we denote by $\pack_\eps \spc{X}$ the maximal number of points in $\spc{X}$ on distance $>\eps$ from each other.
Clearly $\pack_\eps \spc{X}$ is finite for any $\eps>0$.

\begin{thm}{Lemma}\label{lem:pull-back}
Let $\spc{X}$ and $\spc{Y}$ be metric spaces, 
$\spc{X}$ is compact and $f,h\:\spc{X}\to \spc{Y}$ be two continuous maps.

Assume for any $x\in \spc{X}$, $|f(x)-h(x)|<\delta$ then
for any $x,x'\in \spc{X}$ we have
$$\dist_{f,\eps}(x,x')
\le 
\dist_{h,\eps}(x,x')+4\cdot\delta\cdot\pack_\eps \spc{X}.$$

\end{thm}

\parit{Proof.}
Assume $\dist_{h,\eps}(x,x')<\ell$; that is, there is an $\eps$-chain $\{x_i\}_{i=0}^n$ from $x$ to $x'$ such that
$$\sum_{i=1}^n|h(x_{i-1})- h(x_i)|_{\spc{Y}}<\ell.\eqno(*)$$
Since $|h(x_i)-f(x_i)|<\delta$,
$$\dist_{f,\eps}(x,x')\le \sum_{i=1}^n|f(x_{i-1})- f(x_i)|_{\spc{Y}}
<\sum_{i=1}^n|h(x_{i-1})- h(x_i)|_{\spc{Y}}+2\cdot n\cdot \delta$$

Assume $n$ is the smallest number for which there is an $\eps$-chain satisfying $(*)$.
It is enough to show that 
$$n< 2\cdot\pack_\eps \spc{X}.$$

If $n\ge 2\cdot\pack_\eps \spc{X}$, there are $i$ and $j$ such that $j-i>1$ and $|x_i-x_j|\le\eps$.
Remove from this chain all elements $x_k$ with $i<k<j$;
that is, consider new $\eps$-chain 
$$x=x_0,\dots,x_{i-1},x_i,x_j,x_{j+1},\dots,x_n=x'$$
By triangle inequality in $\spc{Y}$, the new chain satisfies $(*)$; 
that is, $n$ is not the smallest number, a contradiction.
\qeds

\begin{thm}{Proposition}\label{prop:diam-preimage} 
Let $\spc{X}$ and $\spc{Y}$ be metric spaces, $\spc{X}$ be compact and $\imath\:\spc{X}\to \spc{Y}$ be an intrinsic isometry.

Then given $\eps>0$ there is $\delta=\delta(\imath,\eps)>0$ such that for any connected set $W\i \spc{X}$  
$$\diam \imath(W)<\delta\ \ \Longrightarrow\ \ \diam W<\eps.$$ 

\end{thm}

\parit{Proof.}
Assume contrary; that is, there is a sequence of connected subsets $W_n\z\i \spc{X}$ such that 
$\diam \imath(W_n)\to 0$ as $n\to\infty$ but $\diam W_n>\eps$.
Thus there are two sequences of points $x_n,x_n'\in W_n$ such that $|x_n -x_n'|\ge\eps$. 
Pass to a subsequence of $n$ so that $W_n\to W$ in Hausdorff sense and $x_n\to x$, $x_n'\to x'$.
We obtain a closed connected subset $W\i \spc{X}$ with two distinct points $x$ and $x'$ such that $\imath(W)=p$ for some $p\in \spc{Y}$. 

Since $W$ is connected, 
for any $\eps>0$ there is an $\eps$-chain $(x_i)_{i=0}^n$ from $x$ to $x'$ such that $\imath(x_i)=p$ for all $i$.
Thus, we have $\dist_{\imath,\eps}(x,x')=0$ for any $\eps>0$;
that is, $\dist_{\imath}(x,x')=0$, a contradiction.
\qeds

\section{The proofs}
\label{proofs}

\parit{Proof of the trivial statement (\ref{top-dim}).}
Given $\eps>0$ choose $\delta=\delta(\imath,\eps)$ as in Proposition~\ref{prop:diam-preimage}. 
Since $\dim \EE^d=d$,
 there is a finite open covering $\{U_i\}_{i=1}^n$ of $\imath(\spc{X})$ with multiplicity $\le d+1$ and
such that $\diam U_i<\delta$ for each $i$.

Consider the covering $\{V_\alpha\}$ of $\spc{X}$ by connected components of $\imath^{-1}(U_i)$ for all $i$.
According to Proposition~\ref{int>int}, $\spc{X}$ is a length space.
In particular,
all sets $V_\alpha$ are open.
According to Proposition~\ref{prop:diam-preimage}, $\diam V_\alpha<\eps$.
Clearly multiplicity of $\{V_\alpha\}$ is at most $d+1$.
Thus, the statement follows.
\qeds

\parit{Proof of ``if'' in the main theorem (\ref{main}).} 
Let $\spc{X}$ be a pro-Euclidean space of rank $\le d$.
Assume $(\spc{P}_n)_{n=0}^\infty$ is a sequence of $d$-dimensional Euclidean polyhedrons and 
$$\phi_{m,n}\: \spc{P}_m\to \spc{P}_n$$
is an inverse system of short maps such that $\spc{X}=\varprojlim \spc{P}_n$.
Let $\psi_n\:\spc{X}\to \spc{P}_n$ be the projections.

According to Theorem~\ref{PL-Nash},
given $\eps_{n+1}>0$ and a piecewise linear intrinsic isometry $\imath_n\:\spc{P}_n\z\to \EE^d$
there is a piecewise linear intrinsic isometry 
$\imath_{n+1}\:\spc{P}_{n+1}\z\to \EE^d$
such that
$$|\imath_{n+1}(x)-\imath_{n}{\circ}\phi_{n+1,n}(x)|<\eps_{n+1}.$$
for any $x\in \spc{P}_n$.
It remains to show that sequence $\eps_{n}$ can be chosen on such a way that $\imath_n{\circ}\psi_n$ converges to an intrinsic isometry $\imath\:\spc{X}\to\EE^d$.

Let us choose 
$\eps_{n+1}>0$ so that 
$$\eps_{n+1}<\tfrac12\min\l\{\eps_n,\delta(\imath_n,\tfrac1n)\r\},$$
 where $\delta(\imath_n,\tfrac1n)$ as in Proposition~\ref{prop:pull-back}.
Clearly, $\sum_i\eps_i<\infty$, thus the the following limit exists
$$\imath=\lim_{n\to\infty} \imath_n\circ\psi_n,\ \ \imath\:\spc{X}\to\EE^d.$$
Obviously, $\imath$ is short.
Further, for any $x\in \spc{X}$,
$$|\imath(x)-\imath_n{\circ}\psi_n(x)|< \sum_{i=n+1}^\infty\eps_i<\delta(\imath_n,\tfrac1n).$$ 
Thus, according to Proposition~\ref{prop:pull-back},
$$\dist_{\imath}(x,x')+\tfrac1n
>
\dist_{\imath_n\circ\psi_n}(x,x')
\ge
|\psi_n(x)-\psi_n(x')|_{\spc{P}_n}.$$
Since $|\psi_n(x)-\psi_n(x')|_{\spc{P}_n}\to |x-x'|_{\spc{X}}$ as $n\to\infty$, the map $\imath\:\spc{X}\to\EE^d$ is an intrinsic isometry.\qeds

\parit{Proof of ``only if'' in the main theorem (\ref{main}).} 
We will give a construction a polyhedron $\spc{P}$ 
associated to an intrinsic isometry $\imath\: \spc{X}\to \EE^d$ and a tiling of $\EE^d$ by coordinate $a$-cubes.
(The space $\spc{P}$ will be glued out of $a$-cubes.)
The construction will be done in such a way that if a tiling $\tau'$ is a subdivision of a tiling $\tau$ then for corresponding polyhedrons $\spc{P}'$ and $\spc{P}$ there will be a natural intrinsic isometry $\spc{P}'\to \spc{P}$. 
Thus we will construct the needed inverse system of polyhedrons out of nested subdivisions of $\EE^d$.

Take sequences $a_n=\tfrac{1}{2^{n}}$ and set $r_n=\tfrac{1}{10}\cdot a_n$.
Fix $n$ for a while and consider tiling of $\EE^d$ by coordinate $a_n$-cubes.
Let us construct a Euclidean polyhedron $\spc{P}_n$ associated to this tiling.

The image $\imath(\spc{X})$ is covered by finite number of such $a_n$-cubes, say $\{\square_{n}^i\}$.
For each $\square_{n}^i$, consider all connected components $\{W^{i j}_n\}$ of 
$$B_{r_n}(\imath^{-1}(\square_{n}^i))\i \spc{X},$$
where $B_r(S)$ denotes $r$-neighborhood of set $S$.

According to Proposition~\ref{int>int},  $\spc{X}$ is a length space.
In particular, each set $W^{i j}_n$ is open and contains a ball of radius $r_n$.
Thus for fixed $i$ the collection of open sets $\{W^{i j}_n\}$ is finite .
Therefore the set of all $\{W^{i j}_n\}$ for all $\{\square^{i}_n\}$ forms a finite open cover of $\spc{X}$.
For each $W^{i j}_n$ make an isometric copy $\square^{i j}_n$ of $\square^{i}_n$ and fix an isometry $\imath^{i j}_n\:\square^{i j}_n\to\square^{i}_n$.
The Euclidean polyhedron $\spc{P}_{n}$, is glued from $\square^{i j}_n$ by the following rule:
glue $\square^{i_1j_1}_n$ to $\square^{i_2j_2}_n$ along $(\imath^{i_2j_2}_n)^{-1}\circ\imath^{i_1j_1}_n$ if and only if $W^{i_1j_1}_n\cap W^{i_2j_2}\not=\emptyset$.
(The map $(\imath^{i_2j_2}_n)^{-1}\circ\imath^{i_1j_1}_n$ sends one of the faces of $\square^{i_1j_1}_n$ isometrically to a  face of $\square^{i_2j_2}_n$.)

The constructed polyhedron $\spc{P}_{n}$ admits a natural piecewise linear intrinsic isometry $\imath_n\:\spc{P}_{n}\to\EE^d$,
defined as 
$\imath_n(x)=\imath^{i j}_n(x)$ if $x\in \square^{i j}_n$.
Further, there is uniquely defined intrinsic isometry $\phi_{m,n}\:\spc{P}_m\to \spc{P}_n$ for $m\ge n$ which satisfies $\imath_m=\imath_n\circ\phi_{m,n}$ and
$$\phi_{m,n}(\square^{i' j'}_m)\i\square^{i j}_n\i \spc{P}_n
\ \ \Rightarrow\ \ 
W^{i' j'}_m\i W^{i j}_n\i \spc{X}.$$
Further, set $\psi_n\:\spc{X}\to \spc{P}_n$ to be intrinsic isometry which uniquely determined by
$\imath_n\circ\psi_{n}=\imath$ and 
$$ \psi_{n}(x)\in\square^{i j}_n\i \spc{P}_n\ \ \Rightarrow\ \ x\in W^{i j}_n\i \spc{X}.$$
Clearly, $\spc{P}_n$ together with $\phi_{m,n}$ form an inverse system and
$\psi_n=\phi_{m,n}\circ\psi_m$ for any pair $m\ge n$.

In order to prove that $\spc{X}=\varprojlim \spc{P}_n$, 
it only remains to show that
$$|x-x'|_{\spc{X}}
\le
\lim_{n\to\infty}|\psi_n(x)-\psi_n(x')|_{\spc{P}_n}\eqno(*)$$
for all $x,x'\in\spc{X}$.

Given a subset $K\i \spc{P}_n$, 
let us denote by $K^*\i \spc{X}$ the union of all $W^{i j}_n\i \spc{X}$ such that $\square^{i j}_n\cap K\not=\emptyset$.
Clearly, if $K$ is connected then so is $K^*$.
More over, $\imath(K^*)\i B_{r_n}(\imath_n(K))$.
Thus, from Proposition~\ref{prop:diam-preimage}, we have that for any $\eps>0$ we can find $\delta>0$ such that 
$$r_n+\diam K<\delta\ \ \Longrightarrow\ \ \diam K^*<\eps \eqno(**)$$

Assume $(*)$ is wrong,
then one can choose $x,x'\in\spc{X}$ and $\eps,\ell>0$ so that
\newcommand*{\threestar}{\mathrel{\vcenter{\offinterlineskip\hbox{$\mkern4.5mu {*}$}\hbox{${*}{*}$}}}}
$$\dist_{\imath,\eps}(x,x')>\ell> |\psi_n(x)-\psi_n(x')|_{\spc{P}_n}\eqno(\threestar)$$
for all $n$.
In particular, for all $n$ there is a path $\gamma_n\:[0,1]\to \spc{P}_n$ from $\psi_n(x)$ to $\psi_n(x')$ with length $<\ell$.
Choose $\delta=\delta(\imath,\eps)$ as in Proposition~\ref{prop:diam-preimage}.
Let $0\z=t_0\z<t_1\z<\dots<t_m=1$ be such that 
\newcommand*{\fourstar}{%
\mathrel{\vcenter{\offinterlineskip
\hbox{${*}{*}$}\hbox{${*}{*}$}}}}
$$\diam\gamma([{t_{i-1}},{t_i}])<\tfrac{\delta}{2}.\eqno(\fourstar)$$ 
Clearly one can assume that 
$m\le2\cdot\lceil\tfrac{\ell}{\delta}\rceil$.
For each $t_i$ choose a point $x_i\z\in \gamma(t_i)^*\i \spc{X}$; clearly 
\newcommand*{\fivestar}{%
\mathrel{\vcenter{\offinterlineskip
\hbox{$\mkern4.5mu {*}{*}$}\hbox{${*}{*}{*}$}}}}
$$|\imath(x_i)-\imath_n{\circ}\gamma(t_i)|_{\EE^d}
<2\cdot a_n.
\eqno(\fivestar)$$
Note that $x_{i-1},x_i\in \gamma([t_{i-1},t_i])^*$.
Thus, $(\fourstar)$ and $(**)$ imply that 
$$|x_{i-1}-x_i|<\diam \gamma_n([t_{i-1},t_i])^*<\eps$$
for all large $n$.
Thus $x_i$ forms an $\eps$-chain from $x$ to $x'$, and $(\fivestar)$ implies
\begin{align*}
 \dist_{\imath,\eps}(x,x')
&\le 
\sum_{i=1}^m|\imath(x_{i-1})-\imath(x_i)|
<\\
&<
\sum_{i=1}^m|\imath_n{\circ}\gamma_n(t_{i-1})-\imath_n{\circ}\gamma_n(t_i)|
+4\cdot a_n\cdot \lceil\tfrac{\ell}{\delta}\rceil
<\\
&<\ell+4\cdot a_n\cdot \lceil\tfrac{\ell}{\delta}\rceil
\end{align*}
which contradicts $(\threestar)$ for large enough $n$.
\qeds

\parit{Remark.} In the constructed inverse system 
$(\phi_{m,n},\spc{P}_n)$, the images of $\phi_{m,n}$ form a $\sqrt{d}\cdot a_n$-net in $\spc{P}_n$.
It follows that the space $\spc{X}$ is isometric to the Gromov--Hausdorff limit of $\spc{P}_n$ (see also Section~\ref{preliminaries}).

\parit{Proof of \ref{minkowski}.}
The inequality $d\le m$ follow from trivial statement (\ref{top-dim}).
In the proof of the second part, we use the following two statements:

\begin{enumerate}
 \item Assume $\imath\:\Omega\i\MM^d \to\EE^m$ is an intrinsic isometry, then it is a Lipschitz map for a Euclidean structure on $\Omega$.
Thus, according to Rademacher's theorem (see \cite[3.1.6]{ref-to-df}) the differential $d_p\imath$ is well defined almost all $p\in \Omega$.
 \item For any curve $\gamma(t)$ with natural parameter in a metric space, we have that for almost all values of parameter $t_0$ we have
$$|\gamma({t_0})-\gamma({t_0+\eps})|=\eps+o(\eps),$$
see  \cite[2.7.5]{BBI}.
\end{enumerate}

Let us denote by $\|{*}\|$ the norm which induces metric on $\MM^d$.
Fix $u$ so that $\|u\|=1$.
Consider pencil of lines of the form $p+u\cdot t$ in $\Omega$.
Two statements above imply that $|d_p\imath(v)|\z\ae\|{v}\|$
Hence we obtain  parallelogram identity 
$$2\cdot\l(\|u\|^2+\|v\|^2\r)=\|u+v\|^2+\|u-v\|^2$$
is satisfied for any two vectors $v$ and $w$.
That is, the norm $\|{*}\|$ is Euclidean.
\qeds

\section{About path isometries}\label{path.isometry}

In this section we will compare the notion of intrinsic isometry defined in Section~\ref{preliminaries} with more common (but less natural) notion of path and weak path isometries.

\begin{thm}{Definition}\label{def:path-iso}
Let $\spc{X}$ and $\spc{Y}$ be two length spaces.
A map $\imath\: \spc{X}\to \spc{Y}$ is called 
\begin{enumerate}
\item\emph{path isometry} if for any path $\gamma\:[0,1]\to \spc{X}$ we have
$$\length \gamma=\length (\imath\circ\gamma).$$
\item\emph{weak path isometry} if for any rectifiable path $\gamma\:[0,1]\to \spc{X}$ we have
$$\length \gamma=\length (\imath\circ\gamma).$$
\end{enumerate}
\end{thm}

As it was noted in Section~\ref{preliminaries}, any intrinsic isometry is a path isometry 
(and therefore, a weak path isometry).
Next we will show that converse does not hold.
Similar counterexamples for weak path isometries are much simpler:
one can take a left-invariant sub-Riemannian metric $d$ on Heisenberg group $H$ then factorizing by center gives an weak path isometry $(H,d)\to\EE^2$ (which is not a path isometry and therefore not an intrinsic isometry).

\begin{thm}{Example}\label{exam:path} There is a compact length space $\spc{X}$ and a path isometry $f\:\spc{X}\z\to \RR$ such that $f^{-1}(0)$ is a  nontrivial connected subset.
Moreover, the Lebesgue covering dimension of $f^{-1}(0)$ can be made arbitrary large.
\end{thm}

In particular, an analog of \ref{top-dim} does not hold for path isometries.

The following construction was suggested by Dmitri Burago;
it is based on two ideas: 
(1) the construction in \cite[3.1]{BIS},
(2) the construction of \emph{pseudo-arc} given in \cite{knaster} (see also the survey \cite{pseudo.arc} and references therein).
In fact, for the first part of theorem $f^{-1}(0)$ will be homeomorphic to a pseudo-arc 
and for the second part $f^{-1}(0)$ will be homeomorphic to a product of pseudo-arcs.

\begin{wrapfigure}{r}{32mm}
\begin{lpic}[t(-5mm),b(3mm),r(0mm),l(0mm)]{pics/pseudoarc(0.2)}
\lbl[b]{140,5;$\II$}
\lbl[l]{16,200;$\JJ$}
\lbl[r]{6,25;$\eps$}
\lbl[b]{85,-17;Graph of an}
\lbl[t]{85,-18;$\eps$-crooked map.}
\end{lpic}
\end{wrapfigure}

Recall that an onto map $h\:\II\to \JJ$ between two real intervals is called \emph{$\eps$-crooked} if for any two values $t_1\z<t_2$ in $\II$ there are values $t_1\z<t'_2\z<t'_1\z<t_2$ such that
$|h(t_i')-h(t_i)|\le\eps$ for $i\in\{1,2\}$.
The existence of $\eps$-crooked map for any given $\eps>0$ is easy to prove by induction on $n=\lceil\tfrac1\eps\cdot{\length\JJ}\rceil$.

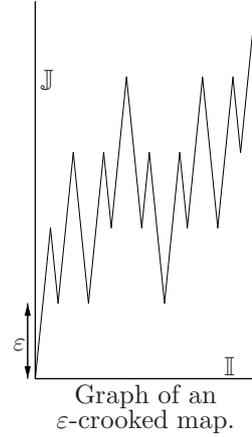
\parit{Proof.}
Start with a metric graph\footnote{That is, locally finite graph with length metric, such that each edge is isometric to a real interval.} $\Gamma$;
consider its completion $\bar\Gamma$.
Set $\grave\Gamma=\bar\Gamma\backslash\Gamma$.
Consider the map $f\:\bar\Gamma\to\RR$, where $f(x)$ is the distance from $x$ to $\grave\Gamma$.
Note that $f$ is a path isometry on $\Gamma$ and $f(\grave\Gamma)=0$.

Let us construct $\Gamma$  such that $\grave\Gamma$ is connected and contains a pair of points that cannot be connected by a curve $\alpha$ such that the curve $f\circ\alpha$ has arbitrary small length. 

Fix a sequence $\eps_n$ that converges to $0$ very fast.
Consider a sequence of real intervals $\JJ_n$ with short $\eps_n$-crooked maps $h_n\:\JJ_{n}\to \JJ_{n-1}$.
We can assume that $\JJ_0=[-1,1]$.

The topological inverse limit $\JJ_\infty=\varprojlim \JJ_n$ is a connected compact space with no nontrivial paths; in fact $\JJ_\infty$ is a pseudo-arc.

For each $n$, choose an $\eps_n$-dense set of \emph{vertexes} in $\JJ_n$.
Connect each vertex $x$ in $\JJ_n$ to the point $h_n(x)$ in  $\JJ_{n-1}$ by an edge of length $\tfrac1{2^n}$.
Set $\Gamma$ to be the obtained graph.

Let us denote by $\Gamma_n$ the finite subgraph of $\Gamma$ formed by $\JJ_0,\JJ_2,\dots,\JJ_{n-1}$ and all the edges between them.
Note that there is a short map $\phi_n\:\bar\Gamma\to\Gamma_{n}$ that is identity on $\Gamma_{n}$ and such that
\[\phi_n|_{\JJ_m}=h_n\circ\dots\circ h_{m}\:\JJ_m\to\JJ_{n-1}\] 
for $m\ge n$.
In particular, $\phi_n$ maps the $\grave\Gamma$ onto $\JJ_{n-1}$.

Assume that the ends of the pseudoarc $\grave\Gamma$ can be connected by a path $\alpha\:[0,1]\to\bar\Gamma$ such that 
\[\length f\circ\alpha<\tfrac1{10}.\]
Without loss of generality, we may assume that $\alpha$ is a simple curve; that is, it has no self-intersections.

Let construct a nested sequence of arcs $\alpha=\alpha_1\supset \alpha_2\supset \dots$ such that $\alpha_n\subset \bar\Gamma\backslash\Gamma_n$ and $\phi_1(\alpha_n)$ contains a segment of length $1$ for each $n$;
in particular diameter of $\alpha_n$ is at least $1$.

Set $c_1=1$; consider the sequence defined by $c_{n+1}=c_n-n\cdot \eps_n$.
Since $\eps_n$ converges to zero very fast, we can assume that $c_n>\tfrac12$ for any $n$.
Therefore it is sufficient to construct arcs $\alpha_1\supset \alpha_2\supset \dots$ such that $\phi_1(\alpha_n)\supset [-c_n,c_n]$ for any $n$.

Clearly we can take $\alpha_1=\alpha$.
Assume the sequence is constructed up to $\alpha_n$, so $\phi_1(\alpha_n)\supset [-c_n,c_n]$.
By the definition of a crooked map,
$\alpha_n$ contains $3$ disjoint subarcs such that $\phi_1$ maps each onto $[-c_n+\eps_n,c_n-\eps_n]$;
each of these arcs contains $3$ disjoint subarcs such that $\phi_1$ maps each  onto $[-c_n+2\cdot\eps_n,c_n-2\cdot \eps_n]$ and so on.
It follows that there are $3^n$ disjoint subarcs $\beta_1,\dots \beta_{3^n}$ such that $\phi_1(\beta_i)\supset [-c_{n+1},c_{n+1}]$.
Note that $\alpha_n$ runs at most $2^n$ times  from $\JJ_n$ to $\JJ_{n+1}$
and each visit to $\JJ_n$ is short.
Therefore $\beta_i\subset \bar\Gamma\backslash\Gamma_{n+1}$ for some $i$; denote this arc by $\alpha_{n+1}$.

Note that the intersection $\bigcap_n\alpha_n$ is an arc in $\grave\Gamma$ of diameter at least $1$, but $\grave\Gamma$ has no nontrivial arcs --- a contradiction.

Let $d$ be the length metric induced by $f$ on $\bar\Gamma$;
that is, the distance $d(x,y)$ is defined as the exact lower bound on the lengths of $f\circ\alpha$ for all paths $\alpha$ from $x$ to $y$.
Note that the map $f\:(\bar\Gamma,d)\to\RR$ is a path isometry and the set $\grave\Gamma$ remains connected in $(\bar\Gamma,d)$; hence the first part follows

\parit{Second part.}
We construct a graph $\Gamma^{(m)}$ to make $\grave\Gamma^{(m)}$ homeomorphic to a product on $m$ copies of $\grave\Gamma$.

We will do the case $m=2$;
the others are analogous.
The set of vertexes of $\Gamma^{(2)}$ is disjoint union $\sqcup_n (\mathop{\rm Vert}\JJ_n\times \mathop{\rm Vert}\JJ_n)$, where $\mathop{\rm Vert}\JJ_n$ denotes the set of vertexes of $\JJ_n$.
We connect two vertexes $(x,y)\in \mathop{\rm Vert}\JJ_n\times \mathop{\rm Vert}\JJ_n$ and $(x',y')\in \mathop{\rm Vert}\JJ_k\times \mathop{\rm Vert}\JJ_k$ if and only if the pairs $(x,x')$  and $(y,y')$ were connected in $\Gamma$; the length of this edge must be maximum of lengths of edges $xx'$ and $y y'$ (we assume that a vertex is connected to it-self by an edge of length $0$).

Clearly, there is a homeomorphism $\grave\Gamma^{(2)}\to\grave\Gamma\times\grave\Gamma$.
Note that there are two short coordinate projections $\varsigma_1,\varsigma_2\:\Gamma^{(2)}\to\Gamma$ which can be extended to the projections $\bar\varsigma_1,\bar\varsigma_2\:\bar\Gamma^{(2)}\to\bar\Gamma$.
Thus for any path $\alpha\:[0,1]\to\bar\Gamma^2$, we have that total length of $\alpha\backslash\grave\Gamma^{(2)}$ is at least as big as its projections.
Hence the second part follows.
\qeds

\section{Comments and open questions}

A length space $\spc{M}$ is called \emph{Minkowski $d$-polyhedron} if there is a finite triangulation of $\spc{M}$ such that each simplex is isometric to a simplex in a Minkowski space.
Correspondingly, a compact metric space $\spc{X}$ is called \emph{pro-Minkowski space} of rank $\le d$ 
if it can be presented as an inverse limit of Minkowski $d$-polyhedrons.

\begin{thm}{Question}
Is it true that any length space with Lebesgue's covering dimension $d$ is a pro-Minkowski space of rank $d$? 
\end{thm}

Or even more specific:

\begin{thm}{Question}\label{mink-disk}
Is it true that any metric space which homeomorphic to a disk is a pro-Minkowski space of rank $2$?
\end{thm}

One can reformulate it philosophically: 
\textit{Is there any essential difference between Finsler metric and general metric on $n$-manifold?}
This question was asked by Dmitri Burago; 
it was also original motivation for this paper
(see also a related example \cite[theorem 1]{BIS}).

If one removes restriction on dimension, 
then the answer to the above question is ``yes''.
Namely, the following exercise can be solved by using Kuratowski embedding $x\mapsto \mathop{\rm dist}_x$.

\begin{thm}{Exercise}
Show that any compact length space is an inverse limit 
of Minkowski polyhedrons $\spc{M}_n$ with $\dim \spc{M}_n\to\infty$.
\end{thm}

\begin{thm}{Question}
Is it true that any path isometry from a closed Euclidean ball to Euclidean space is an intrinsic isometry? 
\end{thm}

The answer is ``yes'' in 2-dimensional case, a proof can be build on the idea suggested by Taras Banakh in \cite{banakh}.

\end{document}

\title{Внутренние изометрии \\
в евклидово пространство}
\author{Антон Петрунин}
\date{}
\maketitle

\begin{abstract}
Я рассматриваю пространства,
допускающие «внутренние изометрии» в $d$-мерное евклидово пространство.
Основной результат состоит в том, что класс таких пространств совпадает с классом обратных пределов $d$-мерных полиэдров. 
\end{abstract}

\section{Введение}

Внутренние изометрии определены в разделе~\ref{ru-preliminaries}, это вариация понятия \emph{изометрии на путях}, т.е. отображений, сохраняющих длины кривых.
Любая внутренняя изометрия является изометрией на путях, обратное, вообще говоря, не верно.
Следующее утверждение является одной из причин, почему я предпочитаю понятие внутренней изометрии.

\begin{thm}{Стартовое предложение}\label{ru-top-dim}
Предположим, что компактное метрическое пространство $\spc{X}$ 
допускает внутреннюю изометрию в $d$-мерное евклидово пространство (далее обозначаемое $\EE^d$).
Тогда 
$\dim \spc{X}\le d$,
где $\dim$ обозначает размерность Лебега.
\end{thm}

Это предложение доказывается в разделе \ref{ru-proofs}.
Пример~\ref{ru-exam:path} показывает,
что для изометрий на путях это утверждение не верно.
Хаусдорфова размерность $\spc{X}$ также не может быть ограничена.
Например, $\RR$-дерево допускает внутреннюю изометрию в $\RR$,
при этом оно содержит компактные подмножества произвольно большой хаусдорфовой размерности.

\medskip

Вот несколько известных результатов о пространствах, допускающих внутренние изометрии в $\EE^d$.

\begin{thm}{Теорема}\label{ru-thm:gromov}
Пусть $\spc{R}$ есть $d$-мерное риманово пространство 
и $f\:\spc{R}\z\to\EE^d$ --- короткое отображение%
\footnote{то есть $1$-липшицево}%
.
Тогда для любого $\eps>0$ 
существует внутренняя изометрия $\imath\:\spc{R}\to\EE^d$
такая,что
$$|f(x)-\imath(x)|_{\EE^d}<\eps$$
для любого $x\in \spc{R}$.

В частности, любое $d$-мерное риманово пространство 
допускает внутреннюю изометрию в $\EE^d$.
\end{thm}

Для изометрий на путях эта теорема была доказана Михаилом Громовым \cite[2.4.11]{gromov-PDE},
но доказательство работает для внутренних изометрий  без каких-либо изменений.
Из этой теоремы следует, что произвольный предел возрастающей последовательности римановых метрик на фиксированном многообразии допускает внутреннюю изометрию в $\EE^d$,
доказывается аналогично условию достаточности в основной теореме (\ref{ru-main}).
В частности, любая субриманова метрика на $d$-мерном многообразии допускает внутреннюю изометрию в $\EE^d$.

\begin{thm}{Теорема}\label{ru-PL-Nash}
Пусть $\spc{P}$ есть полиэдр и $f\: \spc{P}\to\EE^d$ --- короткое отображение.
Тогда для любого $\eps>0$ существует кусочно-линейная внутренняя изометрия $\imath\: \spc{P}\to \EE^d$ 
такая, что
$$|f(x)-\imath(x)|_{\EE^d}<\eps$$
для любой точки $x\in \spc{P}$.
\end{thm}

\begin{thm}{Следствие}\label{ru-zalgaller}
Любой $d$-мерный полиэдр допускает кусочно-линейную внутреннюю изометрию в $\EE^d$.
\end{thm}

Следствие~\ref{ru-zalgaller} было доказано Виктором Залгаллером \cite{ru-zalgaller} для размерностей $\le 4$, 
но небольшая модификация доказательства работает во всех размерностях,
это было показано Светланой Крат \cite{ru-krat}.
Двумерный случай вышеприведённой теоремы был также доказан Крат.
Позже Арсений Акопян \cite{ru-akopjan} обобщил это доказательство на все размерности,
используя кусочно-линейный аналог теоремы Нэша --- Кёйпера.
Эта теорема была доказана Ульрихом Брэмом \cite{ru-brehm},
но его работа на многие годы была оставлена без видимого внимания
и передоказана независимо Арсением Акопяном и Алексеем Тарасовым \cite{ru-akopjan-tarasov}.

\parbf{Необходимое и достаточное условие.} 
Компактное метрическое пространство $\spc{X}$ 
называется \emph{проевклидовым пространством ранга $\z\le d$},
если оно может быть представлено как  \emph{обратный предел}\footnote{определение обратного предела дано в  разделе~\ref{ru-preliminaries}} $\spc{X}=\varprojlim \spc{P}_n$ 
последовательности $d$-мерных полиэдров $\spc{P}_n$.

\begin{thm}{Основная теорема}\label{ru-main}
Компактное метрическое пространство $\spc{X}$ 
допускает внутреннюю изометрию в $\EE^d$ тогда и только тогда, когда $\spc{X}$ является проевклидовым пространством  ранга $\le d$. 
\end{thm}

Формулировка этой теоремы более интересная нежели её доказательство ---  по-моему, это первый случай, когда обратные пределы решают естественную геометрическую задачу.

Из основной теоремы следует, что утверждение теоремы \ref{ru-thm:gromov} (для компактных пространств) эквивалентно 
тому, что любое $d$-мерное риманово пространство является проевклидовым пространством ранга $d$.
Последнее утверждение следует напрямую из упражнения ниже.
В частности, основная теорема даёт альтернативное доказательство теоремы  \ref{ru-thm:gromov} для компактных пространств.

\begin{thm}{Упражнение}
Докажите, что любое компактное риманово пространство допускает липшицеву аппроксимацию полиэдрами.
\end{thm}

\parbf{Не пример.} 
Напомним, что пространство Минковского --- это конечномерное векторное пространство с метрикой, индуцированной некоторой нормой.

\begin{thm}{Предложение}\label{ru-minkowski}
Пусть $\Omega$ есть открытое подмножество $d$-мерного пространства Минковского $\MM^d$.
Предположим, что $\Omega$ допускает внутреннюю\footnote{то же верно и для изометрий на путях} изометрию в $\EE^m$, тогда $d\le m$ и $\MM^d$ изометрично $\EE^d$.
\end{thm}
В частности, условие в \ref{ru-top-dim} на размерность Лебега не достаточно.

\smallskip

{\sloppy

\parbf{Благодарности.}
Автор выражает признательность Арсению Акопяну, Дмитрию Бураго, Сергею Иванову, Александру Лычаку и анонимному референту за полезные письма и беседы.

}

\section{Предварительные замечания}\label{ru-preliminaries}

\parbf{Соглашения.}
Расстояние между точками $x$ и $x'$ метрического пространства $\spc{X}$ будет обозначаться $|x -x'|$ или $|x -x'|_{\spc{X}}$.

Отображение $f\:\spc{X}\to \spc{Y}$ между метрическими пространствами $\spc{X}$ и $\spc{Y}$ называется \emph{коротким},
если  
$$|f(x)-f(x')|_{\spc{Y}}\le |x -x'|_{\spc{X}}$$
для любых $x,x'\in \spc{X}$.

Пространство $\spc{P}$ с внутренней метрикой называется 
\emph{$d$-мерным полиэдром}, если существует конечная триангуляция  $\spc{P}$, каждый симплекс в которой изометричен симплексу в $\EE^d$
.

\parbf{Обратный предел.}
Рассмотрим обратную систему компактных метрических пространств
$(\spc{X}_n)_{n=0}^\infty$ и коротких отображений $\phi_{m,n}:\spc{X}_m\to \spc{X}_n$ при $m\ge n$;
то есть:
\begin{enumerate}
\item $\phi_{m,n}\circ \phi_{k,m}=\phi_{k,n}$ для любой тройки $k\ge m\ge n$ и
\item для любого $n$ отображение $\phi_{n,n}\:\spc{X}_n\to \spc{X}_n$ тождественно.
\end{enumerate}
Компактное метрическое пространство $\spc{X}$ 
является \emph{обратным пределом системы $(\phi_{m,n}, \spc{X}_n)$} (для краткости $\spc{X}=\varprojlim \spc{X}_n$), если его подлежащее множество состоит из всех последовательностей $x_n\in \spc{X}_n$ таких, что $\phi_{m,n}(x_m)\z=x_n$ для всех $m\ge n$,
и для двух таких последовательностей $(x_n)$ и $(x'_n)$ расстояние определяется как 
$$|(x_n)-(x'_n)|_{\spc{X}}=\lim_{n\to\infty}|x_n- x'_n|_{\spc{X}_n}.$$

Пусть $\spc{X}=\varprojlim \spc{X}_n$, тогда отображения $\psi_n\:\spc{X}\to \spc{X}_n$, определяемые как $\psi_n\:(x_i)_{i=0}^\infty\mapsto x_n$ называются \emph{проекциями}.
Таким образом, $\psi_n=\phi_{m,n}\circ\psi_m$ для всех $m\ge n$. 

\parit{Комментарии.}
Вышеприведённое определение эквивалентно обычному обратному пределу в категории, в которой объекты --- компактные метрические пространства,  
а морфизмы --- короткие отображения.

Заметьте, что обратный предел не всегда определён, и
если определён, то результат компактен по определению.
(В принципе, категорию компактных метрических пространств можно расширить так, чтобы предел любой обратной системы был определён.)

Нетрудно видеть, что обратный предел пространств с внутренней метрикой также имеет внутреннюю метрику.

Вообще говоря, обратный предел системы пространств может отличаться от предела по Громову --- Хаусдорфу.
Например, рассмотрим обратную систему $\spc{X}_n=[0,1]$, с отображениями $\phi_{m,n}(x)\equiv 0$. 
Обратный предел этой системы изометричен одноточечному пространству, тогда как предел по Громову --- Хаусдорфу изометричен отрезку $[0,1]$.
Тем не менее нетрудно видеть, что если для любого $\eps>0$ образы $\phi_{m,n}$ образуют $\eps$-сеть в $X_n$ для всех достаточно больших $m$ и $n$ , то $\spc{X}=\varprojlim \spc{X}_n$ изометричен пределу по Громову --- Хаусдорфу.

\parbf{Внутренние изометрии и поднятия метрик.}
Пусть $\spc{X}$ и $\spc{Y}$ суть метрические пространства, 
и $f\:\spc{X}\to \spc{Y}$ есть 
непрерывное отображение.
Для точек $x,x'\in \spc{X}$ последовательность точек $x=x_0,x_1,\dots,x_n=x'$ в $\spc{X}$ называется $\eps$-цепью из $x$ в $x'$, если $|x_{i-1}-x_i|\le\eps$ для всех $i>0$.
Определим
$$\dist_{f,\eps}(x,x')
=
\inf\left\{\sum_{i=1}^n|f(x_{i-1})-f(x_{i})|_{\spc{Y}}\right\}$$
где точная нижняя грань взята по всем $\eps$-цепям $(x_i)_{i=0}^n$ из $x$ в $x'$.

Заметим, что $\dist_{f,\eps}$ есть псевдометрика%
\footnote{то есть $\dist_{f,\eps}$ удовлетворяет неравенству треугольника, она симметрична, неотрицательна и $\dist_{f,\eps}(x,x)=0$, но может случиться, что $\dist_{f,\eps}(x,x')=0$ при $x\not=x'$.}
на $\spc{X}$.
Далее, 
$\dist_{f,\eps}(x,x')$ не возрастает по $\eps$.
Значит, существует (возможно бесконечный) предел
$$\dist_{f}(x,x')=\lim_{\eps\to0+}\dist_{f,\eps}(x,x').$$
Псевдометрика $\dist_f\:\spc{X}\times\spc{X}\to[0,\infty]$ будет называться \emph{поднятием} метрики по $f$.

Отображение $f\:\spc{X}\to \spc{Y}$ называется \emph{внутренней изометрией}, если 
$$|x-x'|_{\spc{X}}=\dist_f(x,x')$$
для любых $x,x'\in \spc{X}$.

Любая внутренняя изометрия является коротким отображением.
Более того, легко видеть, что внутренние изометрии сохраняют длины кривых.
Обратное не верно, см. раздел~\ref{ru-path.isometry}.

\begin{thm}{Предложение}\label{ru-int>int}
Пусть $\spc{X}$ --- компактное (или даже ограниченно компактное\footnote{То есть все замкнутые ограниченные множества в $\spc{X}$ компактны.}) метрическое пространство.
Тогда существование внутренней изометрии $f\:\spc{X}\to \spc{Y}$ влечёт то, что метрика на $\spc{X}$ является внутренней.
\end{thm}

Доказательство этого предложения предоставляется читателю.
Заметим, что в общем случае это не верно.
Рассмотрим две точки, соединённые счётным числом копий единичных отрезков $\II_n$ и одним отрезком длины $\frac12$ с естественной внутренней метрикой.
Выбросим из этого пространства отрезок длины $\frac12$,
метрика на полученном пространстве $\spc{X}$ не является внутренней.
Далее построим отображение $f\:\spc{X}\to\RR$ так, что сужение $f_n=f|_{\II_n}$ является внутренней изометрией, $f_n(0)=0$, $f_n(1)=\tfrac12$ и $f_n(x)$ равномерно сходится к $\tfrac x2$.
Нетрудно видеть, что $f\:\spc{X}\to\RR$ есть внутренняя изометрия.

Следующее предложение аналогично полунепрерывности функционала длины.

\begin{thm}{Предложение}\label{ru-prop:pull-back} 
Пусть $\spc{X}$ и $\spc{Y}$ --- метрические пространства, 
$\spc{X}$ компактно и  непрерывное отображение $f\:\spc{X}\to \spc{Y}$ таково, что 
$$\sup_{x,x'\in\spc{X}}\dist_f(x,x')<\infty.$$ 
Тогда для любого $\eps>0$ существует $\delta=\delta(f,\eps)>0$ такое, что
$$|f(x)-h(x)|_{\spc{Y}}<\delta\ \ \text{для любой точки}\ \ x\in \spc{X}$$
влечёт 
$$\dist_{f}(x,x')<\dist_{h}(x,x')+\eps\ \ \text{для любых}\ \ x,x'\in \spc{X}.$$

\end{thm}

Предложение следует напрямую из леммы~\ref{ru-lem:pull-back}.

Для компактного метрического пространства $\spc{X}$ и $\eps>0$ 
обозначим через $\pack_\eps \spc{X}$ максимальное число точек в $\spc{X}$ на расстоянии $>\eps$ друг от друга.
Очевидно, что $\pack_\eps \spc{X}$ конечно.

\begin{thm}{Лемма}\label{ru-lem:pull-back}
Пусть $\spc{X}$ и $\spc{Y}$ --- метрические пространства, 
$\spc{X}$ компактно и $f,h\:\spc{X}\to \spc{Y}$ --- два
непрерывных отображения.

Предположим, что $|f(x)-h(x)|<\delta$ для любого $x\in \spc{X}$.
Тогда $$\dist_{f,\eps}(x,x')
\le 
\dist_{h,\eps}(x,x')+4\cdot\delta\cdot\pack_\eps \spc{X}.$$
для любых $x,x'\in \spc{X}$.
\end{thm}

\parit{Доказательство.}
Предположим, $\dist_{h,\eps}(x,x')<\ell$, то есть существует $\eps$-цепь $\{x_i\}_{i=0}^n$ из $x$ в $x'$ такая, что
$$\sum_{i=1}^n|h(x_{i-1})- h(x_i)|_{\spc{Y}}<\ell.\eqno(*)$$
Поскольку $|h(x_i)-f(x_i)|<\delta$,
$$\dist_{f,\eps}(x,x')\le \sum_{i=1}^n|f(x_{i-1})- f(x_i)|_{\spc{Y}}
<\sum_{i=1}^n|h(x_{i-1})- h(x_i)|_{\spc{Y}}+2\cdot n\cdot \delta$$

Предположим, $n$ есть минимальное число, для которого существует $\eps$-цепь, удовлетворяющая $(*)$.
Достаточно показать, что 
$$n< 2\cdot\pack_\eps \spc{X}.$$

Если $n\ge 2\cdot\pack_\eps \spc{X}$, то найдутся $i$ и $j$ такие, что $j-i>1$ и $|x_i-x_j|\le\eps$.
Выбросим из цепи все элементы $x_k$ с $i<k<j$,
то есть рассмотрим новую $\eps$-цепь 
$$x=x_0,\dots,x_{i-1},x_i,x_j,x_{j+1},\dots,x_n=x'$$ 
По неравенству треугольника в $\spc{Y}$ новая цепь удовлетворяет $(*)$.
Значит, $n$ не минимально --- противоречие.
\qeds

\begin{thm}{Предложение}\label{ru-prop:diam-preimage} 
Пусть $\spc{X}$ и $\spc{Y}$ --- метрические пространства, 
$\spc{X}$ компактно и $\imath\:\spc{X}\to \spc{Y}$ --- внутренняя изометрия.

Тогда для любого $\eps>0$ существует $\delta=\delta(\imath,\eps)>0$ такое, что для любого связного множества $W\subset \spc{X}$  
$$\diam \imath(W)<\delta\ \ \Longrightarrow\ \ \diam W<\eps.$$ 

\end{thm}

\parit{Доказательство.}
Предположим противное, то есть существует последовательность связных пространств $W_n\subset \spc{X}$ такая, что $\diam \imath(W_n)\to 0$ и $\diam W_n>\eps$.
Значит, существуют две последовательности точек $x_n,x_n'\z\in W_n$ такие, что $|x_n -x_n'|\ge\eps$. 
Поскольку $\spc{X}$ компактно,
можно перейти к подпоследовательности $n$ так, что $x_n\to x$, $x_n'\to x'$ и $W_n\z\to W$ в смысле Хаусдорфа.
Заметим, что связность $W_n$ влечёт связность $W$.
Таким образом, мы получаем замкнутое связное множество $W\subset \spc{X}$ с двумя различными точками $x$ и $x'$ такое, что $\imath$ постоянно на $W$. 

Так, для любого $\eps>0$ существует $\eps$-цепь $(x_i)_{i=0}^n$ из $x$ в $x'$ такая, что $\imath(x_0)=\imath(x_1)=\dots=\imath(x_n)$.
Отсюда $\dist_{\imath,\eps}(x,x')=0$ для любого $\eps>0$.
То есть, $\dist_{\imath}(x,x')=0$ --- противоречие.
\qeds

\section{Доказательства}
\label{ru-proofs}

\parit{Доказательство стартового предложения (\ref{ru-top-dim}).}
Пусть $\eps>0$.
Выберем $\delta\z=\delta(\imath,\eps)$ как в предложении~\ref{ru-prop:diam-preimage}. 
Поскольку $\dim \EE^d=d$, существует конечное открытое покрытие $\{U_i\}$ образа $\imath(\spc{X})$ с кратностью $\le d+1$ и
такое, что $\diam U_i<\delta$ для всех $i$.

Рассмотрим покрытие $\{V_\alpha\}$ пространства $\spc{X}$ всеми связными компонентами всех множеств $\imath^{-1}(U_i)$.
Согласно предложению \ref{ru-int>int} пространство $\spc{X}$ имеет внутреннюю метрику.
В частности,
все множества $V_\alpha$ открыты.
Заметим, что кратность $\{V_\alpha\}$ не превосходит кратности $\{U_i\}$ 
и $\diam V_\alpha<\eps$ (последнее следует из предложения~\ref{ru-prop:diam-preimage}).
\qeds

\parit{Доказательство достаточности условия в \ref{ru-main}.} 
Пусть $\spc{X}$ есть проевклидово пространство ранга $\le d$.
Предположим, $(\spc{P}_n)_{n=0}^\infty$ есть последовательность $d$-мерных полиэдров и  
$\phi_{m,n}\: \spc{P}_m\to \spc{P}_n$
--- обратная система коротких отображений такая, что $\spc{X}=\varprojlim \spc{P}_n$.
Обозначим через $\psi_n\:\spc{X}\to \spc{P}_n$ соответствующие проекции.

Согласно теореме~\ref{ru-PL-Nash}
для любого $\eps_{n+1}>0$ и кусочно-линейной внутренней изометрии $\imath_n\:\spc{P}_n\z\to \EE^d$
существует кусочно-линейная внутренняя изометрия 
$\imath_{n+1}\:\spc{P}_{n+1}\z\to \EE^d$ такая, что неравенство
$$|\imath_{n+1}(x)- \imath_{n}{\circ}\phi_{n+1,n}(x)|<\eps_{n+1}$$
выполняется для любой точки $x\in \spc{P}_n$.
Остаётся показать, что последовательность $\eps_{n}$ может быть выбрана так, 
что $\imath_n{\circ}\psi_n$ сходится к внутренней изометрии $\imath\:\spc{X}\z\to\EE^d$.

Выберем  
$\eps_{n+1}>0$ так, что 
$$\eps_{n+1}<\tfrac12\min\left\{\eps_n,\delta(\imath_n,\tfrac1n)\right\},$$
где $\delta(\imath_n,\tfrac1n)$ как в предложении~\ref{ru-prop:pull-back}.
Ясно, что $\sum_i\eps_i<\infty$. 
Отсюда существует предел
$$\imath=\lim_{n\to\infty} \imath_n\circ\psi_n,\ \ \imath\:\spc{X}\to\EE^d.$$
Очевидно, что $\imath$ --- короткое.
Далее,
$$|\imath(x)- \imath_n{\circ}\psi_n(x)|< \sum_{i=n+1}^\infty\eps_i<\delta(\imath_n,\tfrac1n)$$ 
для любого $x\in \spc{X}$.
Так, согласно предложению~\ref{ru-prop:pull-back},
$$\dist_{\imath}(x,x')+\tfrac1n
>
\dist_{\imath_n\circ\psi_n}(x,x')
\ge
|\psi_n(x)-\psi_n(x')|_{\spc{P}_n}.$$
Поскольку $|\psi_n(x)-\psi_n(x')|_{\spc{P}_n}\to |x-x'|_{\spc{X}}$ при $n\to\infty$, отображение $\imath\:\spc{X}\to\EE^d$ есть внутренняя изометрия.\qeds

\parit{Доказательство необходимости условия в \ref{ru-main}.} 
Мы проведём построение полиэдра $\spc{P}$ 
по внутренней изометрии $\imath\: \spc{X}\to \EE^d$ и замощению $\EE^d$ координатными $a$-кубами,
то есть, кубами со стороной $a$.
(Полиэдр $\spc{P}$ будет склеен из $a$-кубов.)
Построение будет обладать следующим свойством:
если замощение $\tau'$ является подразбиением замощения $\tau$, то для соответствующих полиэдров $\spc{P}'$ и $\spc{P}$ существует естественная внутренняя изометрия $\spc{P}'\to \spc{P}$. 
Таким образом, мы сможем построить необходимую обратную систему полиэдров по последовательности подразбиений одного замощения $\EE^d$.

Пусть $a_n=\tfrac{1}{2^{n}}$ и $r_n=\tfrac{1}{10}\cdot a_n$.
Зафиксируем $n$ на некоторое время.
Рассмотрим замощение $\EE^d$ координатными $a_n$-кубами.

Образ $\imath(\spc{X})$ покрывается конечным семейством $a_n$-кубов $\{\square_{n}^i\}$ из нашего замощения.
Для каждого $\square_{n}^i$ рассмотрим все связные компоненты $\{W^{i j}_n\}$ множества
$B_{r_n}(\imath^{-1}(\square_{n}^i))\subset \spc{X}$,
где $B_r(S)$ обозначает $r$-окрестность множества $S$.

Согласно предложению \ref{ru-int>int}, пространство $\spc{X}$ имеет внутреннюю метрику.
В частности, каждое множество $W^{i j}_n$ открыто и содержит шар радиуса $r_n$.
Отсюда, для фиксированного $i$ 
семейство $\{W^{i j}_n\}$ конечно.
Таким образом, все $W^{i j}_n$ образуют конечное открытое покрытие $\spc{X}$.
Для каждого $W^{i j}_n$ заготовим копию $\square^{i j}_n$ куба $\square^{i}_n$ с изометрией $\imath^{i j}_n\:\square^{i j}_n\to\square^{i}_n$.
Полиэдр $\spc{P}_{n}$ клеится из $\square^{i j}_n$ по следующему правилу:
куб $\square^{i_1j_1}_n$ склеен с $\square^{i_2j_2}_n$ по $(\imath^{i_2j_2}_n)^{-1}\circ\imath^{i_1j_1}_n$ тогда и только тогда, когда $W^{i_1j_1}_n\cap W^{i_2j_2}\not=\emptyset$.
(Заметьте, что $(\imath^{i_2j_2}_n)^{-1}\circ\imath^{i_1j_1}_n$ изометрично отображает одну грань $\square^{i_1j_1}_n$ на грань $\square^{i_2j_2}_n$.)

{\sloppy 

Построенный полиэдр $\spc{P}_{n}$ допускает естественную кусочно-линейную внутреннюю изометрию $\imath_n\:\spc{P}_{n}\to\EE^d$,
определяемую как
$\imath_n(x)=\imath^{i j}_n(x)$ при $x\in \square^{i j}_n$.
Далее, существует однозначно определённая внутренняя изометрия  $\phi_{m,n}\:\spc{P}_m\to \spc{P}_n$ при $m\ge n$, которая удовлетворяет условиям $\imath_m=\imath_n\circ\phi_{m,n}$ и такая, что
$$\phi_{m,n}(\square^{i' j'}_m)\subset\square^{i j}_n\subset \spc{P}_n
\ \ \Rightarrow\ \ 
W^{i' j'}_m\subset W^{i j}_n\subset \spc{X}.$$
Проекции $\psi_n\:\spc{X}\z\to \spc{P}_n$ однозначно определяются условиями
$\imath_n\circ\psi_{n}=\imath$ и
$$ \psi_{n}(x)\in\square^{i j}_n\subset \spc{P}_n\ \ \Rightarrow\ \ x\in W^{i j}_n\subset \spc{X}.$$
Очевидно, $\spc{P}_n$ вместе с $\phi_{m,n}$ образуют обратную систему, и
$\psi_n=\phi_{m,n}\circ\psi_m$ при всех $m\ge n$.

}

Чтобы доказать $\spc{X}=\varprojlim \spc{P}_n$, 
достаточно проверить, что
$$|x-x'|_{\spc{X}}
\le
\lim_{n\to\infty}|\psi_n(x)-\psi_n(x')|\eqno(*)$$
для всех $x,x'\in\spc{X}$.

Для множества $K\subset \spc{P}_n$ 
обозначим через $K^*\subset \spc{X}$ объединение всех $W^{i j}_n\subset \spc{X}$ таких, что $\square^{i j}_n\cap K\not=\emptyset$.
Заметим, что если $K$ связно, то связно и $K^*$.
Более того, $\imath(K^*)\subset B_{r_n}(\imath_n(K))$.
Так, из предложения~\ref{ru-prop:diam-preimage} 
получаем, что для любого $\eps>0$  существует $\delta>0$ такое, что 
$$r_n+\diam K<\delta\ \ \Longrightarrow\ \ \diam K^*<\eps \eqno(**)$$

Предположим $(*)$ не верно,
тогда можно выбрать $x,x'\in\spc{X}$ и $\eps,\ell>0$ такие, что
\newcommand*{\threestar}{\mathrel{\vcenter{\offinterlineskip\hbox{$\mkern4.5mu {*}$}\hbox{${*}{*}$}}}}
$$\dist_{\imath,\eps}(x,x')>\ell> |\psi_n(x)-\psi_n(x')|_{\spc{P}_n}\eqno(\threestar)$$
при всех $n$.
В частности, для любого $n$ существует путь $\gamma_n\:[0,1]\to \spc{P}_n$ из $\psi_n(x)$ в $\psi_n(x')$ длины $<\ell$.
Выберем $\delta=\delta(\imath,\eps)$ как в предложении~\ref{ru-prop:diam-preimage}.
Пусть значения $0\z=t_0\z<t_1\z<\dots<t_m=1$ таковы, что 
\newcommand*{\fourstar}{%
\mathrel{\vcenter{\offinterlineskip
\hbox{${*}{*}$}\hbox{${*}{*}$}}}}
$$\diam\gamma([{t_{i-1}},{t_i}])<\tfrac{\delta}{2}.\eqno(\fourstar)$$ 
Можно также предположить, что
$m\le2\cdot\lceil\tfrac{\ell}{\delta}\rceil$.
Для каждого $t_i$ выберем точку $x_i\z\in \gamma(t_i)^*\subset \spc{X}$.
Ясно, что 
\newcommand*{\fivestar}{%
\mathrel{\vcenter{\offinterlineskip
\hbox{$\mkern4.5mu {*}{*}$}\hbox{${*}{*}{*}$}}}}
$$|\imath(x_i)- \imath_n{\circ}\gamma(t_i)|_{\EE^d}
<2\cdot a_n.
\eqno(\fivestar)$$
Отметим, что $x_{i-1},x_i\in \gamma([t_{i-1},t_i])^*$.
Значит, $(\fourstar)$ и $(**)$ влекут 
$$|x_{i-1}-x_i|<\diam \gamma_n([t_{i-1},t_i])^*<\eps$$
для всех достаточно больших $n$.
Отсюда, точки $x_i$ образуют $\eps$-цепь из $x$ в $x'$, и $(\fivestar)$ влечёт
\begin{align*}
 \dist_{\imath,\eps}(x,x')
&\le 
\sum_{i=1}^m|\imath(x_{i-1})-\imath(x_i)|
<\\
&<
\sum_{i=1}^m|\imath_n{\circ}\gamma_n(t_{i-1})-\imath_n{\circ}\gamma_n(t_i)|
+4\cdot a_n\cdot \lceil\tfrac{\ell}{\delta}\rceil
<\\
&<\ell+4\cdot a_n\cdot \lceil\tfrac{\ell}{\delta}\rceil,
\end{align*}
что противоречит $(\threestar)$ для достаточно больших $n$.
\qeds

\parit{Замечание.} В построенной обратной системе 
$(\phi_{m,n},\spc{P}_n)$, образы $\phi_{m,n}$ образуют $\sqrt{d}\cdot a_n$-сеть в $\spc{P}_n$.
Из этого следует, что пространство $\spc{X}$ изометрично пределу $\spc{P}_n$ по Громову --- Хаусдорфу (см. также раздел \ref{ru-preliminaries}).

\parit{Доказательство предложения~\ref{ru-minkowski}.}
Неравенство $d\le m$ следует из стартового предложения \ref{ru-top-dim}.
Для доказательства второй части нам потребуются два утверждения:

\begin{enumerate}
\item Предположим, $\imath\:\Omega\subset\MM^d \to\EE^m$ есть внутренняя изометрия, 
тогда $\imath$ --- липшицево для евклидовой метрики на $\Omega$.
По теореме Радемахера (см. \cite[3.1.6]{ref-to-df}) дифференциал $d_p\imath$ определён для почти всех $p\in \Omega$.
\item Пусть $\gamma(t)$ --- кривая с натуральным параметром в метрическом пространстве.
Тогда 
$$|\gamma({t_0})-\gamma({t_0+\eps})|=\eps+o(\eps) $$
для почти всех значений параметра $t_0$,
см.  \cite[2.7.5]{BBI}.
\end{enumerate}
Обозначим через $\|{*}\|$ норму, которая индуцирует метрику на $\MM^d$.
Зафиксируем вектор $u$ такой, что $\|u\|=1$.
Рассмотрим пучок кривых вида $p+u\cdot t$ в $\Omega$.
Применив оба вышеприведённых утверждения, получаем,
что $|d_p\imath(v)|\z\ae\|{v}\|$.
Отсюда следует тождество параллелограмма
$$2\cdot\left(\|v\|^2+\|w\|^2\right)=\|v+w\|^2+\|v-w\|^2$$
для любых векторов $v$ и $w$.
То есть норма $\|{*}\|$ евклидова.
\qeds

\section{Об изометриях на путях}\label{ru-path.isometry}

В этом разделе сравнивается понятие внутренней изометрии с более распространённым (и менее естественным) понятием изометрии на путях и слабой изометрии на путях.

\begin{thm}{Определение}\label{ru-def:path-iso}
Пусть $\spc{X}$ и $\spc{Y}$ суть пространства с внутренней метрикой.
Отображение $\imath\: \spc{X}\to \spc{Y}$ называется 
\begin{enumerate}
\item\emph{изометрией на путях}, если для любой кривой $\gamma$ в $\spc{X}$ выполняется 
$$\length \gamma=\length \imath\circ\gamma.$$
\item\emph{слабой изометрией на путях}, если для любой спрямляемой кривой $\gamma$ в $\spc{X}$ выполняется
$$\length \gamma=\length \imath\circ\gamma.$$
\end{enumerate}
\end{thm}

Как было отмечено в разделе~\ref{ru-preliminaries}, любая внутренняя изометрия является изометрией на путях 
(в частности, слабой изометрией на путях).
Далее мы покажем, что обратное не верно.
Для слабой изометрии на путях есть очень простой пример:
рассмотрим левоинвариантную субриманову метрику $d$ на группе Гейзенберга $H$.
Тогда факторизация по центру есть слабая изометрия на путях $(H,d)\to\EE^2$, которая не является изометрией на путях и уж тем более  внутренней изометрией.

{

\begin{wrapfigure}{r}{32mm}
\begin{lpic}[t(-9mm),b(4mm),r(0mm),l(0mm)]{pics/pseudoarc(0.2)}
\lbl[b]{140,5;$\II$}
\lbl[l]{16,200;$\JJ$}
\lbl[r]{6,25;$\eps$}
\lbl[b]{85,-16;График $\eps$-скючен-}
\lbl[t]{85,-17;ного отображения.}
\end{lpic}
\end{wrapfigure}

\begin{thm}{Пример}\label{ru-exam:path} Существует компактное пространство с внутренней метрикой $\spc{X}$ и изометрия на путях $f\:\spc{X}\to \RR$ такая, что $f^{-1}(0)$ --- связное неодноточечное множество.

Более того, в таком примере размерность Лебега у $f^{-1}(0)$ может быть сделана произвольно большой.
\end{thm}

В частности, аналог \ref{ru-top-dim} не выполняется для изометрий на путях.

Нижеприведённое построение пространства $\spc{X}$ подсказал мне Дмитрий Бураго;
оно использует две идеи: 
(1) построение в \cite[3.1]{BIS},
(2) построение \emph{псевдодуги Кнастера} в \cite{ru-knaster} (см. также обзор \cite{ru-pseudo.arc}).
В этом построении $f^{-1}(0)$ гомеоморфно произведению псевдодуг.

}

Напомним, что сюръекция $h\:\II\to \JJ$ между вещественными интервалами называется \emph{$\eps$-скрюченной} если для любых $t_1\z<t_2$ в $\II$ найдутся $t_1\z<t'_2\z<t'_1\z<t_2$ такие, что
$|h(t_i')-h(t_i)|\le\eps$ для $i\in\{1,2\}$.
Существование $\eps$-скрюченного отображения для данного $\eps>0$ легко доказывается индукцией по $n\z=\lceil\tfrac1\eps\cdot{\length\JJ}\rceil$.

\parit{Доказательство.}
Начнём с метрического графа\footnote{То есть локально конечного графа с внутренней метрикой такой, что каждое ребро изометрично вещественному интервалу.} $\Gamma$ и рассмотрим его пополнение $\bar\Gamma$.
Пусть $\grave\Gamma=\bar\Gamma\backslash\Gamma$.
Рассмотрим отображение $f\:\bar\Gamma\to\RR$, где $f(x)$ есть расстояние от $x$ до $\grave\Gamma$.
Заметим, что $f$ является изометрией на путях в $\Gamma$ и $f(\grave\Gamma)=0$.

Далее мы построим граф $\Gamma$ такой, что $\grave\Gamma$ связанно и содержит пару точек, которые нельзя соединить путём $\alpha$ таким, что путь $f\circ\alpha$ имеет произвольно малую длину. 

Выберем последовательность положительных чисел $\eps_n$, которая быстро  сходится к нулю.
Рассмотрим последовательность интервалов $\JJ_n$ с короткими $\eps_n$-скрюченными отображениями $h_n\:\JJ_{n}\to \JJ_{n-1}$.
Можно предположить, что $\JJ_0=[-1,1]$.

Топологический обратный предел $\JJ_\infty=\varprojlim \JJ_n$ есть связанное компактное пространство без нетривиальных путей; на самом деле $\JJ_\infty$ есть псевдодуга.

Для каждого $n$, выберем $\eps_n$-плотное множество \emph{вершин} в $\JJ_n$.
Соединим каждую вершину $x$ в $\JJ_n$ с $h_n(x)$ в $\JJ_{n-1}$ ребром длины $\tfrac1{2^n}$.
Полученный граф обозначим $\Gamma$.

Обозначим через $\Gamma_n$ конечный подграф в $\Gamma$ образованный $\JJ_0,\JJ_2,\dots,\JJ_{n-1}$ и всеми рёбрами между ними.
Заметим, что существует короткое отображение $\phi_n\:\bar\Gamma\to\Gamma_{n}$ тождественное на $\Gamma_{n}$ и такое, что
\[\phi_n|_{\JJ_m}=h_n\circ\dots\circ h_{m}\:\JJ_m\to\JJ_{n-1}\] 
при $m\ge n$.
В частности, $\phi_n$ отображает $\grave\Gamma$ на $\JJ_{n-1}$.

Предположим, что концы псевдодуги $\grave\Gamma$ можно соединить путём $\alpha$ таким, что
\[\length f\circ\alpha<\tfrac1{10}.\]
Не умоляя общности, можно предположить, что $\alpha$ простая кривая; то есть, не имеет самопересечений.

Построим убывающую последовательность дуг $\alpha=\alpha_1\supset \alpha_2\supset \dots$ такую, что $\alpha_n\subset \bar\Gamma\backslash\Gamma_n$ и $\phi_1(\alpha_n)$ содержит интервал длины $1$ для каждого $n$;
в частности диаметр $\alpha_n$ хотя бы $1$.

Положим $c_1=1$; 
рассмотрим последовательность определённую рекурсивным соотношением $c_{n+1}=c_n-n\cdot \eps_n$.
Поскольку $\eps_n$ быстро сходится к нулю, можно предположить, что $c_n>\tfrac12$ для любого $n$.
Таким образом, достаточно построить дуги $\alpha_1\supset \alpha_2\supset \dots$ такие, что $\phi_1(\alpha_n)\supset [-c_n,c_n]$ для любого $n$.

Можно взять $\alpha_1=\alpha$.
Предположим последовательность построена до $\alpha_n$, в частности $\phi_1(\alpha_n)\supset [-c_n,c_n]$.
Из определения скрюченного отображения следут, что
$\alpha_n$ содержит $3$ непересекающихся дуги таких, что $\phi_1$ отображает каждую на $[-c_n+\eps_n,c_n-\eps_n]$;
каждая из этих дуг содержит $3$ дуги таких, что $\phi_1$ отображает каждую на  $[-c_n+2\cdot\eps_n,c_n-2\cdot \eps_n]$ и так далее.
Отсюда следует, что существуют $3^n$ дуги $\beta_1,\dots \beta_{3^n}$ таких, что $\phi_1(\beta_i)\supset [-c_{n+1},c_{n+1}]$ для каждого $i$.
Заметим, что $\alpha_n$ может пройти не более $2^n$ раз из $\JJ_n$ в $\JJ_{n+1}$
и каждый визит в $\JJ_n$ не длинный.
Таким образом $\beta_i\subset \bar\Gamma\backslash\Gamma_{n+1}$ для некоторого $i$; возьмём эту дугу за $\alpha_{n+1}$.

Заметим, что пересечение $\bigcap_n\alpha_n$ есть дуга в $\grave\Gamma$ диаметра не меньше $1$, но $\grave\Gamma$ не содержит нетривиальных дуг --- противоречие.

Рассмотрим метрику $d$, такую, что $d(x,y)$ есть точная нижняя грань длин $f\circ\alpha$ для всех путей $\alpha$ из $x$ в $y$.
Заметим, что отображение $f\:(\bar\Gamma,d)\to\RR$ есть изометрия на путях и множество $\grave\Gamma$ остаётся связным в $(\bar\Gamma,d)$; отсюда следует первая часть.

\parit{Вторая часть.}
Мы построим граф $\Gamma^{(m)}$ так, что $\grave\Gamma^{(m)}$  гомеоморфен произведению $m$ копий $\grave\Gamma$.

Для простоты рассмотрим только случай $m=2$ --- остальные случаи аналогичны.
Множество вершин $\Gamma^{(2)}$ есть $\sqcup_n (\mathop{\rm Vert}\JJ_n\times \mathop{\rm Vert}\JJ_n)$, 
где $\mathop{\rm Vert}\JJ_n$ обозначает множество вершин в $\JJ_n$.
Соединим две вершины $(x,y)\in \mathop{\rm Vert}\JJ_n\times \mathop{\rm Vert}\JJ_n$ и $(x',y')\in \mathop{\rm Vert}\JJ_k\times \mathop{\rm Vert}\JJ_k$, если пары $(x,x')$ и $(y,y')$ соединены в $\Gamma$, и припишем этому ребру длину, равную максимуму длин рёбер   $x x'$ и $y y'$ (мы считаем, что каждая вершина подсоединена к себе ребром длины $0$).

По построению $\grave\Gamma^{(2)}$ гомеоморфно произведению псевдодуг.
Отметим, что естественные проекции $\varsigma_1,\varsigma_2\:\Gamma^{(2)}\to\Gamma$ являются короткими;
их можно продолжить до коротких проекций $\bar\varsigma_1,\bar\varsigma_2\:\bar\Gamma^{(2)}\to\bar\Gamma$.
Значит, для любого пути $\alpha\:[0,1]\to\bar\Gamma^2$ 
общая длина $\alpha\backslash\grave\Gamma$ не может быть меньше их проекций.
Отсюда следует вторая часть.
\qeds

\section{Замечания и открытые вопросы}

Пространство $\spc{M}$ с внутренней метрикой называется \emph{$d$-мерным полиэдром Минковского},
если существует конечная триангуляция $\spc{M}$ такая, что каждый симплекс изометричен симплексу в пространстве Минковского.
Соответственно, компактное метрическое пространство $\spc{X}$
 называется \emph{проминковским пространством} ранга $\le d$, 
если оно является обратным пределом $d$-мерных полиэдров Минковского.

\begin{thm}{Вопрос}
Верно ли, что любое компактное пространство с внутренней метрикой и размерностью Лебега $d$ является проминковским пространством ранга $d$? 
\end{thm}

Более конкретный вопрос:

\begin{thm}{Вопрос}\label{ru-mink-disk}
Верно ли, что любое метрическое пространство гомеоморфное диску есть проминковское пространство ранга 2?
\end{thm}

Этот вопрос можно переформулировать по-философски: 
\textit{Есть ли существенная разница между пространствами всех внутренних метрик на двумерном многообразии и всех финслеровых метрик на нём же?}
Этот вопрос, поставленный Дмитрием Бураго, положил начало настоящей статье (см. также \cite[теорема 1]{BIS}).

Заметим, что вложение Куратовского $x\mapsto \mathop{\rm dist}_x$ даёт решение следующего упражнения.

\begin{thm}{Упражнение}
Покажите, что любое компактное пространство с внутренней метрикой является обратным пределом последовательности полиэдров Минковского $\spc{M}_n$ с $\dim \spc{M}_n\to\infty$.
\end{thm}

\begin{thm}{Вопрос}
Верно ли, что любая изометрия на путях из замкнутого евклидова шара в евклидово пространство является внутренней изометрией? 
\end{thm}

В двумерном случае ответ «да»;
доказательство можно построить на идее Тараса Банаха \cite{ru-banakh}.